\documentclass[a4paper, 12pt]{article}
\usepackage{adjustbox}
\usepackage{graphicx}
\usepackage{amsfonts}
\usepackage{amsbsy}
\usepackage{amssymb}
\usepackage[fleqn]{amsmath}
\usepackage{multicol}
\usepackage{longtable,ltxtable,booktabs}
\usepackage{blkarray}
\usepackage{afterpage}
\usepackage{float}
\usepackage{caption}
\usepackage{subcaption}
\usepackage{multirow}
\usepackage{pdflscape}
\usepackage{latexsym}
\usepackage{rotating}
\usepackage{enumerate}

\usepackage{mathtools}

\usepackage{epstopdf}
\usepackage{caption}
\usepackage{xcolor}
\usepackage[nodisplayskipstretch]{setspace}
\usepackage{amsmath}
\usepackage{longtable}
\usepackage[margin=.8 in]{geometry}

\usepackage{verbatim}

\usepackage{mathrsfs}
\usepackage{amsthm}
\usepackage{setspace}
 \usepackage{pdflscape}
\usepackage{colortbl}
\usepackage{latexsym}
\usepackage{amsfonts}
\usepackage{epsfig}
\usepackage{hyperref}
\usepackage{graphics}
\usepackage{float}
\theoremstyle{definition}
\newtheorem{thm}{Theorem}[section]

\newtheorem{defn}{Definition}[section]

\newtheorem{lemma}{Lemma}[section]

\newtheorem{rem}{Remark}[section]

\newtheorem{cor}{Corollary}[section]

\theoremstyle{definition}

\title{ Cut edges and Central vertices of zero divisor graph of the ring of integers modulo n}

\date{{}}\author{Nabajit Talukdar\thanks{Corresponding author.Email address : ntalukdar2000@yahoo.co.in}\\
\small Department of Mathematics\\
	\small Cotton University\\
    \small Guwahati-781001, India\\
	\small and\\
	Niranjan Bora\thanks{ E-mail address : niranjanbora@dibru.ac.in}\\
	\small Department of Mathematics\\
	\small Dibrugarh University\\
	\small Dibrugarh-786004, India}
\begin{document}
\maketitle

\begin{abstract}

The zero divisor graph of a commutative ring $R$ with unity is a
graph whose vertices are the nonzero zero-divisors of the ring, with
two distinct vertices being adjacent if their product is zero. This graph
is denoted by $\Gamma(R)$. In this article we determine the cut-edges
and central vertices in the graph $\Gamma(\mathbb{Z}_{n})$.

 \end{abstract}

  $\mathbf{2010\;Mathematics\;Subject\;Classification:}\;05C25$
	
	$\mathbf{Keywords\;and\;Phrases:}$ Zero-divisor graph, Cut edge, central Vertex

\section{Introduction}

By a graph we mean a finite, undirected, connected graph without loops or multiple edges. We denote by \mbox{$G=(V, E)$} a graph having the set of vertices $V$ and the set of edges $E$. For any vertex $u\in V$, $N(u)$ denotes the set of neighbours of $u$. Terms not defined here are used in the sense of Harary\cite{harary1972graph}. For the results in number theory we follow \cite{niven1991introduction} . All rings in this paper are commutative with unity. The ring of integers modulo $n$ will be denoted by $\mathbb{Z}_{n}$. Following \cite{dummit2004abstract}, for any ring $R$ and
  $x\in R$, the annihilator of $x$, denoted $ann(x)$, is the set $\{y\in R: xy=0\}$. Recall that an element $a$ of a ring $R$ is said to be a zero-divisor if there exists a non-zero element $b$ of $R$ such that $ab = 0$. Let $Z(R)$ denote the set of zero-divisors of a commutative ring $R$. Zero-divisor graphs were introduced by Beck\cite{istvan1988coloring} and have been studied by many authors (e.g. \cite{mulay2002cycles}, \cite{levy2002zero}, \cite{anderson1999zero}).

First we state some results related to number theory.
\begin{lemma}
\label{lem_num_1}
\cite{niven1991introduction}
If $a,k\in \mathbb{N}$, $\text{gcd}(a,a+k)|k$
\end{lemma}

\begin{lemma}
\label{lem_equiv_cong}
\cite{niven1991introduction} For $a, b,m,n \in \mathbb{N}$, 
$an\equiv bn(\text{mod}\  m)$ iff $a\equiv b(\text{mod} \frac{m}{\text{gcd}(m,n)})$
\end{lemma}

\begin{lemma}
\cite{niven1991introduction} For $a, b,m \in \mathbb{N}$, 
The equation $ax+b \equiv 0(\text{mod}\  m)$ has a solution iff $d | b$, where $d=\text{gcd}(a,n)$.
\end{lemma}

\begin{rem}
The neighbours of any vertex $a$ in $\Gamma (Z_{n})$ are the non-zero solutions of the equation $ax \equiv 0(\text{mod}\ n)$. Thus, if $\text{gcd}(a,n)=d$, then the degree of $a$ is $d-1$. Since, $xa \equiv 0(\text{mod}\ n)$ implies that $x(n-a) \equiv 0(\text{mod}\ n)$, any vertex which is adjacent to $a$ is also adjacent to $n-a$.
\end{rem}

Now we state some definitions and results related to a graph.

\begin{defn}
\cite{harary1972graph}
A cut edge of a graph is an edge whose deletion increases the number of
components.
\end{defn}

\begin{lemma}
\label{lem_cut_edge}
\cite{harary1972graph}
    An edge is a cut edge if and only if it belongs to no cycle.
\end{lemma}

For any two vertices $u$ and $v$ of a graph $G$, if there is a path connecting $u$ and $v$, then the distance from $u$ to $v$, written $d(u, v)$, is the length of the shortest path joining $u$ and $v$. 
\begin{defn}
\cite{harary1972graph}
The eccentricity of a vertex $u$  of a graph $G$ is
$\text{max}\{d(u, v) : v \in V (G)\}$. 
\end{defn}

\begin{defn}
\cite{harary1972graph}
A vertex of a graph with minimum eccentricity is called a central vertex of the graph. The set of all the central vertices of a graph $G$ is called the centre of $G$.
\end{defn}

\begin{defn}
\cite{harary1972graph}
 The diameter of a graph $G$ is defined as $\text{diam}(G) = \text{max}\{d(u, v) : u, v \in V(G)\}$.
\end{defn}

It has been found that the zero divisor graph of a commutative ring is connected. 
\begin{thm}
   \cite{anderson1999zero} 
   Let $R$ be a commutative ring. Then $\Gamma(R)$  is connected and $\text{diam}\ \Gamma(R)\leq 3$.
\end{thm}

In this article we prove the following results.
\begin{thm}
  $(a,b) \in E(\Gamma (\mathbb{Z}_{n}))$ is a cut-edge if and only if $\text{gcd}(a,n)=2$ and $2b=n$.
  \end{thm}

\begin{thm}
For any positive integer $n$, the set of central vertices of the graph $\Gamma (\mathbb{Z}_{n})$  is $\cup_{p} \text{ann}(p)\setminus \{0\}$, where the union is taken over all the primes $p$ which divides $n$.
\end{thm}

\section{Cut edges in $\Gamma(\mathbb{Z}_{n})$}
First we determine which edges of the graph $\Gamma (Z_{n})$ can not be a cut edge.
\begin{lemma}
\label{cut_edge_gcd}
Let $(a,b) \in E(\Gamma (\mathbb{Z}_{n}))$. If $\text{gcd}(a,n) \geq 3$ and $\text{gcd}(b,n) \geq 3$, then $(a,b)$  can not be a cut-edge in $\Gamma (Z_{n})$.
\begin{proof}
Since $(a,b) \in E(\Gamma (Z_{n}))$, $ab\equiv 0(\text{mod}\ n)$. Let $d=\text{gcd}(a,n)$ and $d' =\text{gcd}(b,n)$. By Lemma \ref{lem_equiv_cong}
we get that $a\equiv 0(\text{mod}\ \frac{n}{\text{gcd}(b,n)})$. This gives that 
$ a\equiv 0 (mod\ \frac{n}{d'}) $ and consequently we get that
 $a=\frac{n}{d'}k' $ for some $k' \in Z$. Since $a < n$, $k' < d'$.
  Similarly we get that $b \equiv 0 (\text{mod}\ \frac{n}{d})$ and this in turn gives that $b=\frac{n}{d}k $  for some $k\in Z$ and $k<d$.\\

 The following equation has $d$  solutions including $0$ and $b$.

 \begin{equation}
ax \equiv 0(mod\ n) \label{eq(1)}
\end{equation}
All the solutions of equation (~\ref{eq(1)})are members of the set $A=\{b, b+ \frac{n}{d},b+ 2\frac{n}{d}, \ldots , b+ (d-1)\frac{n}{d}\} $ .\\

Similarly the solutions of the equation
\begin{equation}
by \equiv 0(mod\ n) \label{eq(2)}
\end{equation}
 are members of the set  $B=\{a, a+\frac{n}{d'}, a+2\frac{n}{d'}, \ldots , a+ (d'-1)\frac{n}{d'}\}$.\\
 
 Since $b + (d-1) \frac{n}{d} = \frac{n}{d}(k+d-1) < \frac{n}{d}(2d-1) <2n $ and $a + (d'-1) \frac{n}{d'} = \frac{n}{d'}(k'+d'-1) < \frac{n}{d'}(2d'-1) <2n $, it follows that for every $x\in A\cup B$, $x<2n$. Since $0$ is solution of each of equation (~\ref{eq(1)}) and equation (~\ref{eq(2)}), we have that $n \in A$ and  $n \in B$. Also, $n-b $ is a solution of equation (~\ref{eq(1)}) and $n-a $ is a solution of equation (~\ref{eq(2)}). Hence $n-b \in A $ and $n-a \in B$. 
 Suppose, $n-a\neq n$. We get that $\frac{n}{d}$ and $n-\frac{n}{d}$ are members of $A$ and hence either $b\neq \frac{n}{d} $ or $b \neq n-\frac{n}{d}$. It is to be noted that all the members of $A$(or of $B$), which are distinct from $n$ are neighbours of $a$(or of $b$) in $\Gamma (Z_{n})$. By lemma ~\ref{lem_num_1}, we get $d=(a, a+n-a)|n-a$. Hence, $(n-a)(n - \frac{n}{d})= (n-a)n - n\frac{n-a}{d} \equiv 0(\text{mod}\ n)$ and consequently $n-a$ is adjacent to  $n - \frac{n}{d}$ in $\Gamma (Z_{n})$. Again $(n-a)$ is adjacent to $\frac{n}{d}$ in $\Gamma (Z_{n})$. Thus $n-a ,b, a,n-\frac{n}{d}, n-a$ or $n-a ,b,a,\frac{n}{d}, n-a$ is a cycle according as $b\neq n-\frac{n}{d}$ or $b\neq \frac{n}{d}$. By Lemma ~\ref{lem_cut_edge}, $(a,b)$ is not a cut-edge in $\Gamma (Z_{n})$.
\end{proof}
\end{lemma}

 \begin{thm}
  $(a,b) \in E(\Gamma (\mathbb{Z}_{n}))$ is a cut-edge if and only if $\text{gcd}(a,n)=2$, $\text{gcd}(b,n)\geq 3$ and $2b=n$.
  \begin{proof}
  Suppose $(a,b)$ is a cut edge in $\Gamma (Z_{n})$. Then by Lemma ~\ref{cut_edge_gcd} it follows that either $\text{gcd}(a,n)=2$ or $\text{gcd}(b,n)=2$. If both $\text{gcd}(a,n)=2$ and $\text{gcd}(b,n)=2$, then the vertices $a, b, \frac{n}{2}$ form a cycle in $\Gamma (Z_{n})$. Suppose $\text{gcd}(a,n)=2$ and $\text{gcd}(b,n)\geq 3$. The neighbours of $a$ in $\Gamma (Z_{n})$ are non-zero solutions of the equation $ax \equiv 0(mod\ n)$. Since $\text{gcd}(a,n)=2$, this equation will have only one non zero solution. But $n-b$ and $b$ are both non zero solutions of $ax \equiv 0(mod\ n)$. So, it follows that $2b=n$.\\
  \indent
  Conversely, suppose that $\text{gcd}(a,n)=2$, $\text{gcd}(b,n)\geq 3$  and $2b=n$. Then $b$ is the only neighbour of $2$ in $\Gamma (\mathbb{Z}_{n})$. So removal of the edge  $(a,b)$ increases the number of components in $\Gamma (\mathbb{Z}_{n})$. Hence  $(a,b)$ is a cut-edge.

  \end{proof}
 \end{thm}

\begin{cor}
 For odd value of $n$, $\Gamma (\mathbb{Z}_{n})$ has no cut-edge.
\end{cor}
\




\section{Central vertices in $\Gamma(\mathbb{Z}_{n}$)}
\begin{lemma}
\label{lemma_central_p_k}
If $p$ be a prime and $k>1$ be a positive integer, then the central vertices of $\Gamma(\mathbb{Z}_{p^{k}})$ are $p^{k-1}, 2p^{k-1}, \ldots , (p-1)p^{k-1}$.
\begin{proof}
We get that $V(\Gamma(\mathbb{Z}_{p^{k}}))=\{
p, 2p, \ldots , (p-1)p, p^{2}, 2p^{2}, \ldots , (p-1)p^{2}, \ldots, p^{k-1}, 2p^{k-1}, \ldots , (p-1)p^{k-1}\}$. Each of the vertices in the set $\{p^{k-1}, 2p^{k-1}, \ldots , (p-1)p^{k-1}\}$ is adjacent to all other vertices in $V(\Gamma(\mathbb{Z}_{p^{k}}))$ and hence each of the vertices in the set $\{p^{k-1}, 2p^{k-1}, \ldots , (p-1)p^{k-1}\}$ has eccentricity $1$ . Again for each of the vertices in the set $V(\Gamma(\mathbb{Z}_{p^{k}})) \setminus \{p^{k-1}, 2p^{k-1}, \ldots , (p-1)p^{k-1}\}$, there is vertex which is not adjacent to it. Thus we get that the central vertices of $\Gamma(\mathbb{Z}_{p^{k}})$ are $\{p^{k-1}, 2p^{k-1}, \ldots , (p-1)p^{k-1}\}$.
\end{proof}
\end{lemma}

\begin{lemma}
\label{lemma_central_pq}
Let $n$ be a positive integer such that $n =pq$, where $p$ and $q$  are two distict primes primes. Then the central vertices of $\Gamma(\mathbb{Z}_{n})$ are $p, 2p, \ldots , (q-1)p, q, 2q, \ldots , (p-1)q$.

\begin{proof}
We get that $\Gamma(\mathbb{Z}_{n})$ is complete bi-partite graph where the set of vertices can be partitioned into the sets $\{p, 2p, \ldots , (q-1)p\}$ and $\{q, 2q, \ldots , (p-1)q\}$. Hence it follows that the central vertices of $\Gamma(\mathbb{Z}_{n})$ are $p, 2p, \ldots , (q-1)p, q, 2q, \ldots , (p-1)q$.

\end{proof}
\end{lemma}

\begin{lemma}
 If $p$ and $q$ be two distinct primes dividing $n$ and $n >pq$ then $d(p,q)=3$ in $\Gamma (\mathbb{Z}_{n})$.
 \begin{proof}
 Since $n > pq$, we get that $d(p,q)>1$. Let $d(p,q)=2$. Now the neighbours of $p$ are $\frac{n}{p}$, $2\frac{n}{p}$,\ldots, $(p-1)\frac{n}{p}$ and $q$ must be adjacent to one of these neighbours. Therefore, $qi\frac{n}{p}\equiv 0(\text{mod}n)$ for some $i \in \{1, 2, \ldots, p-1\}$. Then $p|qi$, which is not possible. Hence $d(p,q)=3$ .
 \end{proof}
\end{lemma}

 \begin{cor}
 If   $p$ and $q$ be distinct primes dividing $n$ and $n >pq$, $d(\Gamma (\mathbb{Z}_{n}))=3$.
\end{cor}

\begin{thm}
\label{thm_central_n}
Let $n$ be a positive integer which is divisible by at least two primes $p$ and $q$ and $n >pq$. A vertex $a$  in $\Gamma (\mathbb{Z}_{n})$ is a central vertex if and only if  $a \in \text{ann}(r)$ for some prime $r$ which divides $n$.
\begin{proof}
Let $a\in V(\Gamma(\mathbb{Z}_{n}))$ be such that $a\in ann(r)$ for some prime $r$ which divides $n$.Then  $a\in \{\frac{n}{r},2\frac{n}{r},\ldots, (r-1)\frac{n}{r}\}$ and consequently $\frac{n}{r}|a$. Let $b$ be any other vertex in $\Gamma(\mathbb{Z}_{n})$. We show that either $d(a,b)=1$ or $d(a,b)=2$.
If $p|b$, then $ab\equiv0(\text{mod}\ n)$ and hence $d(a,b)=1$. Let $b$ be not divisible by $p$. Let $b_{0}$ be the smallest positive integer such that $\text{gcd}(n,b_{0})=\text{gcd}(n,b)$. Then $p$ can't divide $b_{o}$ otherwise $\frac{b_{o}}{p}$ would be a positive integer smaller than $b_{o}$  such that $\text{gcd}(n,\frac{b_{o}}{p})=\text{gcd}(n,b)$. Let  $q(\neq p)$ be prime such that $q|b_{0}$. Now, $q|n$,otherwise $\frac{b_{o}}{q}$ would be a positive integer smaller than $b_{o}$  such that $\text{gcd}(n,\frac{b_{o}}{q})=\text{gcd}(n,b_{o})=\text{gcd}(n,b)$. Since $\text{gcd}(n,\frac{b_{o}}{p})=\text{gcd}(n,b)$ and $q|n$, we get that $q|b$. So, $b$ is a neighbour of $\frac{n}{q}$. Since $p,q|n$, $\frac{n}{q}\frac{n}{p}=\frac{n}{pq}n\equiv0(\text{mod}\ n)$. Therefore, $\frac{n}{q}\in N(b)$ is adjacent to $\frac{n}{p}$ and consequently to $a$. Thus $b, \frac{n}{q}, a$ is path and hence $d(a,b)=2$. Thus we get that $a\in V(\Gamma(\mathbb{Z}_{n}))$ is a central vertex.\\

Conversely let $a$  in $\Gamma (\mathbb{Z}_{n})$ be a central vertex. Then for any $b \in V(\Gamma (\mathbb{Z}_{n}))$, $d(a,b)=1$ or $d(a,b)=2$. By Lemma ~\ref{lem_equiv_cong}, we get that the neighbours of $a$ in $\Gamma (\mathbb{Z}_{n})$ are $\frac{n}{d}$, $2\frac{n}{d}$, \ldots, $(d-1)\frac{n}{d}$, where $d=\text{gcd}(a,n)$. We claim: $\frac{n}{d}$ is a prime. Let $r$ be a prime dividing $\frac{n}{d}$. Then either $d(r,a )=1$ or $d(r,a )=2$.
We consider the following cases.
\begin{description}
    \item[Case I:]  $d(r,a )=2$.\\

Then $r$ is adjacent to $i\frac{n}{d}$ for some $i \in \{1, 2, \ldots, d-1\}$.Then since both $r$ and $a$ are neighbours of $i\frac{n}{d}$, $gcd(a,r)>1$ and hence $r|a$.Now, $gcd(a,n)=d$ and $r|a$, $r|n$ imply that $r|d$. Since $r|\frac{n}{d}$ and $r|d$, $n=kr^{2}$ for some $k \in \mathbb{N}$. We claim:  $n$ is not divisible by any prime other than $r$. Since $r|a$,$a$ is adjacent to every neighbour of $r$ and hence $N(r) \subseteq N(a)$. If possible, let $s$ be a prime other than $r$ such that $s|k$. If $d(a,s)=1$, then $s$ is a neighbour of $a$ and hence the only possibility is that $\frac{n}{d}=s$. But this can't happen since $r|\frac{n}{d}$. So, let $d(a,s)=2$. Then $s$ is adjacent to a vertex $N(a) \setminus N(r)$. So, we get that $s|a$. Hence $N(s) \subseteq N(a)$. So, $\{\frac{n}{s},2\frac{n}{s}, \ldots, (s-1)\frac{n}{s}\}\subseteq \{\frac{n}{d},2\frac{n}{d},\ldots, (d-1)\frac{n}{d}\}$.If $s|\frac{n}{d}$, then $s|u $ for all $u\in N(a)$ and in particular $s|\frac{n}{r}$. Then $r, \frac{n}{r}, s$ is a path, which contradicts that $d(r,s)=3$. So, $s$ does not divide $\frac{n}{d}$. Since, every neighbour of $s$ is a multiple of $\frac{n}{d}$ and $s$ does not divide $\frac{n}{d}$, it follows that  $N(s)=\{s\frac{n}{d}, 2s\frac{n}{d}, \ldots, (s-1)\frac{n}{d}\}$ and consequently $\{\frac{n}{s}, 2s\frac{n}{s}, \ldots, (s-1)\frac{n}{s}\}=\{s\frac{n}{d}, 2s\frac{n}{d}, \ldots, (s-1)\frac{n}{d}\}$. This gives that $\frac{n}{s}=s\frac{n}{d}$ and hence $d=s^{2}$. This is a contradiction to the fact that $r|d$. So, we get that $n$ is not divisible by any prime other than $r$. Thus, $n=r^{k}$ for some $k\geq 2$. This contradicts our assumption that there are at least two distinct primes dividing $n$.\\ Hence this case does not arise.

    \item [Case II:]  $d(r,a )=1$.\\
    Now $r$  is a neighbour of $a$ in $\Gamma(\mathbb{Z}_{n})$ and hence $r\in \{\frac{n}{d}, 2\frac{n}{d}, \ldots, (d-1)\frac{n}{d}\}$. This gives that $r=\frac{n}{d}$. Thus $n=dr$ and hence $a\in \text{ann}(r)$.
\end{description}

\end{proof}
\end{thm}

Combining Lemma \ref{lemma_central_p_k}, Lemma \ref{lemma_central_pq} and Theorem \ref{thm_central_n} we get the following result.

\begin{thm}
For any positive integer $n$, the set of central vertices of the graph $\Gamma (\mathbb{Z}_{n})$  is $\cup_{p} \text{ann}(p)\setminus \{0\}$, where the union is taken over all the primes $p$ which divides $n$.
\end{thm}


\begin{thebibliography}{9999}

\bibitem{anderson1999zero} David F. Anderson and Philip S. Livingston , {\it The Zero-Divisor Graph of a Commutative Ring} , {\sf Journal of Algebra}, {\bf 217}:  4.34-447, 1999.

\bibitem{istvan1988coloring}
Istvan Beck, {\it Coloring of Commutative Ring}, {\sf Journal of Algebra}, {\bf 116}:208-226, 1988.


\bibitem{dummit2004abstract}
David S. Dummit and Richard M. Foote, {\it Abstract Algebra}, {\sf John Wiely \& Sons}, 2004.

\bibitem{harary1972graph} F. Harary, {\it Graph Theory}, {\sf  Addition-Wesley, Reading}, 1972.

\bibitem{levy2002zero} Ron Levy and Jay Shapiro, {\it The zero-divisor graph of von Neumann regular rings} , {\sf Communications in Algebra}, {\bf 30(2)}:  3745-750, 2002.



\bibitem{mulay2002cycles}Shashikant B Mulay, {\it Cycles and symmetries of zero-divisors} , {\sf Communications in Algebra}, {\bf 30(7)}:  3533-3558, 2002.



\bibitem{niven1991introduction} I. Niven , H.S. Zuckerman , H.L. Montgomery, {\it }
David S. Dummit and Richard M. Foote, {\it An introduction to the theory of numbers}, {\sf John Wiely \& Sons}, 1991.


\end{thebibliography}
\end{document}